\documentclass[12pt]{amsart}

\theoremstyle{plain}
\newtheorem{thm}          {Theorem}       [section]

\newtheorem{prop} [thm]   {Proposition}

\theoremstyle{definition}
\newtheorem{defn} [thm]   {Definition}

\newtheorem{rk}   [thm]   {Remark}

\theoremstyle{remark}

\renewcommand{\thenotn}{}

\renewcommand\theequation{\arabic{section}.\arabic{equation}}

\newcommand{\bc}{Bochner}
\newcommand{\bct}{Bochner technique}
\newcommand{\bw}{between}
\newcommand{\cc}{\mbox{$connected component$}}
\newcommand{\conf}{conformal foliation}
\newcommand{\dmen}{dimensional}
\newcommand{\h}{harmonic}
\newcommand{\hc}{horizontally conformal}
\newcommand{\hwc}{horizontally (weakly) conformal}
\newcommand{\hm}{harmonic morphism}
\newcommand{\hms}{harmonic morphisms}
\newcommand{\mph}{morphisms of $p$-harmonic functions}
\newcommand{\rman}{Riemannian manifold}
\newcommand{\scon}{simply-connected}
\newcommand{\SF}{Seifert fibre}
\newcommand{\nscon}{non-simply-connected}
\newcommand{\tg}{totally geodesic}
\newcommand{\wf}{Weitzenb\"{o}ck formula}

\newcommand{\lap}{\text{\( \Delta \)}}
\newcommand{\gam}{\text{\( \gamma \)}}
\newcommand{\G}{\text{\( \Gamma \)}}
\newcommand{\Gp}{\text{\( {\Gamma}^{\prime} \)}}
\newcommand{\lm}{\text{\( \lambda \)}}
\newcommand{\m}{\text{\( \mu \)}}
\newcommand{\n}{\text{\( \nabla \)}}
\newcommand{\p}{\text{\( \phi \)}}
\newcommand{\s}{\text{\( \sigma \)}}
\newcommand{\ca}{\text{\( \mathcal A \)}}
\newcommand{\cf}{\text{\( \mathcal F \)}}
\newcommand{\ch}{\text{\( \mathcal H \)}}
\newcommand{\cv}{\text{\( \mathcal V \)}}

\newcommand{\GG}{\mbox{\boldmath \(\mathrm G\)}}
\newcommand{\spn}{\mbox{\rm span}\,} %can't use \span which is reserved
\newcommand{\la}{\left\langle}
\newcommand{\ra}{\right\rangle}
\renewcommand{\angle}{\mbox{\rm angle}} %latex \angle gives nothing

\newcommand{\spf}{\text{\( {\mathbb E}^3 \)}}
\newcommand{\qspf}{\text{\( {\mathbb E}^{3}/\G \)}}

\newcommand{\lo}{\text{\( \lm_1 \)}}     
\newcommand{\lt}{\text{\( \lm_2 \)}}
\newcommand{\li}{\text{\( \lm_i \)}}   
\newcommand{\lj}{\text{\( \lm_j \)}}   
\newcommand{\lk}{\text{\( \lm_k \)}}   

\newcommand{\Xo}{\text{\( X_1 \)}}  
\newcommand{\Xt}{\text{\( X_2 \)}}  
\newcommand{\Ui}{\text{\( U_i \)}}  
\newcommand{\Uj}{\text{\( U_j \)}}

\newcommand{\vol}{\text{\( {{\upsilon}^M} \)}}

\newcommand{\fn}[2]{\text{\( {\p}{\colon}{#1}{\to}{#2} \)}}

\newcommand{\inm}[2]{\text{\( {\boldsymbol\langle}{#1},{#2}{\boldsymbol\rangle}^{M} \)}}
\newcommand{\inn}[2]{\text{\( {\boldsymbol\langle}{#1},{#2}{\boldsymbol\rangle}^{N} \)}}

\newcommand{\sderm}[2]{\text{\( {\n \! _{#1} ^M}\! {#2} \)}}

\newcommand{\derm}[3]{\text{\( \inm{{\n \! _{#1} ^M}\! {#2}}{#3} \)}}

\newcommand{\dderm}[4]{\text{\( \inm{{\n \! _{#1} ^M}\! {\n \! _{#2} ^M}
\! {#3}}{#4} \)}}

\newcommand{\bderm}[4]{\text{\( \inm{{\n \! _{#1} ^M}\! {#2}}
{{\n \! _{#3} ^M}\! {#4}} \)}}

\newcommand{\sdern}{\text{\( {\n \! _{\Ui} ^{\phi ^{\! -1}\! TN}}
\! {d\p{\cdot}X_{k}} \)}}

\newcommand{\dern}{\text{\( \inn{{\n \! _{\Ui} ^{\phi ^{\! -1}\! TN}}
\! {d\p{\cdot}\Xo}}{d\p{\cdot}\Xt} \)}}

\newcommand{\dernjj}{\text{\( \inn{{\n \! _{\Ui} ^{\phi ^{\! -1}\! TN}}
\! {d\p{\cdot}X_{k}}}{d\p{\cdot}X_{k}} \)}}

\newcommand{\e}[2]{\text{\( {({e}_{s})}_{#1}^{#2} \)}}
\newcommand{\es}{\text {\( e_s \)}}
\newcommand{\et}{\text {\( e_t \)}}
\newcommand{\el}{\text {\( e_l \)}}

\newcommand{\ep}[2]{\text{\( {({e}_{s}^{\prime})}_{#1}^{#2} \)}}
\newcommand{\eps}{\text {\( e_s^{\prime} \)}}
\newcommand{\ept}{\text {\( e_t^{\prime} \)}}
\newcommand{\epl}{\text {\( e_l^{\prime} \)}}

\newcommand{\con}[1]{\text{\( \n^{#1} \)}}
\newcommand{\curv}[1]{\text{\( R^{#1} \)}}

\newcommand{\hor}{\text {\( T_x^{H}M \)}}
\newcommand{\ver}{\text {\( T_x^{V}M \)}}

\newcommand{\Rm}{\text {\( \mathbf{R^M} \)}}
\newcommand{\Rn}{\text {\( \mathbf{R^N} \)}}

\newcommand{\ric}[1]{\text{\( \mathbf{{Ricci}^{#1}} \)}}

\newcommand{\riem}[1]{\text{\( \mathbf{{Riem}^{#1}} \)}}

\newcommand{\scalm}{\text{\( \mathbf{{{Scal}^{M}{\mid}_{H}}} \)}}

\newcommand{\scaln}{\text{\( \mathbf{{Scal}^{N}} \)}}

\newcommand{\sm}[2]{\text{\( \sum_{#1}^{#2} \)}}

\newcommand{\ssm}{\text{\( {\displaystyle \sum_{i=1}^{m}} \)}}

\newcommand{\dsm}{\text{\( {\ssm}\sm{j=1,j\neq i}{m} \)}}

\begin{document}
\baselineskip 18pt

%%%%%%%%%%%%%%%%%%%%%%%%%%%%%%%%%%%%%%%%%%%%%%%%%%%%%%%%%%%
%%								
%%
%%		TITLE ETC
%%
%%
%%
%%%%%%%%%%%%%%%%%%%%%%%%%%%%%%%%%%%%%%%%%%%%%%%%%%%%%%%%%%%

\title[Harmonic morphisms from three-dimensional space forms]
{Harmonic morphisms from three-dimensional   
Euclidean and spherical space forms}
 
\author{M.T.~MUSTAFA*}
\author{J.C.~WOOD}



\address{(Until May~1997);\ \  
         Math. Section,
         I.C.T.P., 
         P.O.~Box 586, 
         34100 Trieste, 
         Italy}

\email{mustafa@ictp.trieste.it}

\address{Department of Pure Mathematics,
         University of Leeds,
         Leeds LS2 9JT,
         UK}
         
\email{j.c.wood@leeds.ac.uk}

\thanks{* Supported by the Government of Pakistan and the International 
        Centre for Theoretical Physics, Trieste}

%\date{\today}

%\maketitle

%%%%%%%%%%%%%%%%%%%%%%%%%%%%%%%%%%%%%%%%%%%%%%%%%%%%%%%%%%%
%%								
%%
%%		ABSTRACT
%%
%%
%%
%%%%%%%%%%%%%%%%%%%%%%%%%%%%%%%%%%%%%%%%%%%%%%%%%%%%%%%%%%%

\begin{abstract}
This paper gives a description 
of  all  \hms\ from a 
three-\dmen\ {\it \nscon} 
Euclidean and spherical space form to a surface, by extending 
the work of Baird-Wood \cite{BaiWoo88,BaiWoo91} who 
dealt with the simply-connected case; namely we show that 
any such harmonic morphism is the composition of a ``standard'' 
\hm\ and a weakly conformal map. To complete the description we list  
the space forms and the standard \hm s on them.
\end{abstract}

\maketitle

%%%%%%%%%%%%%%%%%%%%%%%%%%%%%%%%%%%%%%%%%%%%%%%%%%%%%%%%%%%
%%								
%%
%%		FIRST SECTION
%%		
%%		INTRODUCTION
%%
%%%%%%%%%%%%%%%%%%%%%%%%%%%%%%%%%%%%%%%%%%%%%%%%%%%%%%%%%%%

\section{Introduction}
\label{sec:first}

\setcounter{equation}{0}

A smooth map \fn{M}{N} \bw\ \rman s is called a {\it \hm}\ if it preserves 
germs of \h\ functions, i.e., if $f$ is a real-valued \h\ function on 
an open set $V\subseteq N$ then the composition $f\circ \p$ is \h\ on 
${\p}^{-1}(V)\subseteq ~M$.

In \cite{BaiWoo88,BaiWoo91} P.~Baird and the second author studied 
\hms\ from a three-\dmen\ \scon\ space form to a surface and 
obtained a complete local and global classification of them. 
The classification of three-\dmen\ Euclidean and spherical space forms 
is well-known cf. \cite[Chapters~3 and 7]{wol66}. This motivates us 
to study the global classification of \hms\ from three-\dmen\ 
{\it \nscon}\ Euclidean and spherical 
space forms to a surface. In this paper, we obtain a 
 description of {\it all} \hms\ 
from any three-\dmen\ 
 Euclidean and spherical space form to a surface, namely that
any such harmonic morphism is the composition of a 
standard \hm\ (see Remark~\ref{rem:std}) and a weakly conformal map. 

Section~\ref{sec:intro hm} gives some basic facts on harmonic morphisms 
and 
conformal foliations. 
Section~\ref{sec:nscon} contains the main results of the paper.  
Section~\ref{sec:standard} lists the standard \hm s from 
all the spherical and orientable Euclidean space forms. We shall 
make use of standard properties of orbifolds and Seifert fibre spaces; 
the reader 
is referred to \cite{scot83,thur} for these. 

%%%%%%%%%%%%%%%%%%%%%%%%%%%%%%%%%%%%%%%%%%%%%%%%%%%%%%%%%%
%%		SECTION 2
%%        Background material on harmonic morphisms
%%
%%%%%%%%%%%%%%%%%%%%%%%%%%%%%%%%%%%%%%%%%%%%%%%%%%%%%%%%%%
\section{Background material on harmonic morphisms}
\label{sec:intro hm}
\subsection{Harmonic morphisms} \qquad
%\label{subsec:hmorphs}

Let $(M^m,\inm{\cdot}{\cdot})$ and $(N^n,\inn{\cdot}{\cdot})$ be 
smooth $(C^\infty)$ \rman s of dimensions $m,n$ respectively.
\begin{defn}
A smooth map \fn{M^m}{N^n}  is called a 
{\it \hm}\ if,  for every real-valued 
function $f$ which is \h\ on an open subset $V$ of $N$ with $\p^{-1}(V)$ 
non-empty, $f\circ \p$ is a 
real-valued \h\ function on $\p^{-1}(V)\subset M$.
\end{defn}

For a smooth map \fn{M^m}{N^n}, let $C_{\p}
= \{ x\in M\mid {\rm rank}~d\p_x <n \}$ be its 
{\it critical set}. The points of the set $M\setminus C_{\p}$ are 
called {\it regular points}.  For each $x\in M\setminus C_{\p}$, the  
{\it vertical space} $T^{V}_{x}M$ at $x$ is defined by 
$T^{V}_{x}M = {\rm Ker}~d\p_x$. The 
{\it horizontal space} \hor\ at $x$ is given by the 
orthogonal complement of \ver\ in $T_{x}M$ so that 
$T_{x}M = \ver \oplus \hor$.

\begin{defn}
A smooth map \fn{(M^m,\inm{\cdot}{\cdot})}{(N^n,\inn{\cdot}{\cdot})} is called 
{\it horizontally 
 (weakly) conformal}\/ if $d\p = 0$ on 
$C_{\p}$ and the restriction of \p\ to $M\setminus C_{\p}$ is 
a conformal submersion, that is, for each 
$x\in M\setminus C_{\p}$, the differential $d\p_x :\hor \to T_{\p(x)}N$ is conformal 
and surjective. This means that there exists a function 
$\lm:M\setminus C_{\p}\to {\mathbb R}^+$ such that
$$
  \inn{d\p(X)}{d\p(Y)} = \lm^2 \inm{X}{Y} \quad \forall X, Y \in T^{H}M.
$$
\end{defn}

By setting $\lm=0$ on $C_{\p}$, we can extend 
$\lm:M\to {\mathbb R}^{+}_{0}$ to a continuous function on $M$ such that 
${{\lm}^2}$ is smooth, in fact 
${{\lm}^2} = {\|d\phi\|^2}/n$. The function $\lm:M\to {\mathbb R}^{+}_{0}$ 
is called the {\it dilation} of the map \p.

The characterization obtained by 
B.~Fuglede \cite{Fug78} and T.~Ishihara \cite{Ish79}, states that 
{\it \hm s are precisely the \h\ maps which are \hwc}.

When the codomain is two-dimensional, \hms\ \fn{M^m}{N^2} have 
special features. An important property is 
{\it conformal invariance}, explained as follows.
\begin{prop}
Let $N_{1}^{2}$, $N_{2}^{2}$ be surfaces, i.e. 2-dimensional Riemannian manifolds. Let 
\fn{M^m}{N_{1}^{2}} be a \hm\ and 
$\psi :N_{1}^{2}\to N_{2}^{2}$ a weakly conformal map. Then 
$\psi \circ \p :M^m \to N_{2}^{2}$ is a \hm. 
\end{prop}
Thus the notion of a \hm\ to a surface $(N^2,h)$ 
depends only on the conformal equivalence class of $h$. In particular, 
the 
concept of a \hm\ to a Riemann surface is well-defined.

The reader is referred to \cite{Fug78,Bai83,Woo86A}, 
 for further basic properties of \hm s.

\subsection{Harmonic morphisms and conformal foliations}
\label{subsec:h&conf} \qquad

The fibres of a submersion \fn{M^m}{N^n} define a 
foliation on $M$ whose leaves are the connected components of the 
fibres and any foliation is given locally this way. A foliation 
on $(M^m,g)$ is called  {\it conformal} (respectively {\it Riemannian}) 
if it is given locally by conformal (respectively Riemannian) 
submersions from open subsets of $M^m$. For alternative 
definitions and more information see \cite[Sec.~3]{Woo86A}.

The relationship \bw\ \hms\ and \conf s is the following: 
By \cite{Fug78,Ish79} the fibres of 
a submersive \hm\ \fn{M^m}{N^n} define a \conf\ on $M$. In case $N$ 
is a surface,  a 
submersion \fn{M^m}{N^2} is a \hm\ if and only if it is 
\hc\ with minimal fibres \cite{BaiEel81,GUD}. (In fact, this remains true
for any non-constant map, see \cite{Woo96}). When $M$ is a 
3-dimensional manifold, 
the conformal foliation defined by a \hm\ can be extended over the critical 
points of \p\ and we have 
\begin{prop}\cite{BaiWoo92A}
\label{prop:ind conf}
Let \fn{M^3}{N^2} be a non-constant \hm. Then the fibres of \p\ form a 
conformal foliation \cf\ by geodesics of $M^3$.
\end{prop}

It follows that, {\it locally}, $\phi$ is a submersion to an open subset
of $\mathbb C$ followed by a weakly conformal map (cf. \cite{BaiWoo91}).

\subsection{Harmonic morphisms from three-\dmen\ \scon\ 
 space forms to a surface}\label{sec:531} \qquad

A complete classification of \hms\ from a three-\dmen\ {\it \scon}\ 
space form was obtained by P.~Baird and the second author in 
\cite{BaiWoo88,BaiWoo91}. To state their 
 results we consider the following {\it standard examples of 
\hms\ from three-\dmen\ \scon\ space forms}: 
\begin{enumerate}
  \item {\it Orthogonal projection} $\pi_1 :{\mathbb R}^3 \to {\mathbb R}^2$ 
defined by $\pi_1 (x_1,x_2,x_3) = (x_1,x_2)$.
  \item The {\it Hopf map} $\pi_2 :{\mathbb S}^3 \to {\mathbb S}^2$. 
Identifying ${\mathbb S}^2$ with the one-point compactification 
${\mathbb C}\cup \infty$ of the complex numbers via stereographic 
projection and taking ${\mathbb S}^3$ as 
$
 {\mathbb S}^3 = \left\{
   (z_1,z_2)\in{\mathbb C}^2 : \| z_1 \|^2 + \| z_2 \|^2 = 1
             \right\}, 
$ 
 then the Hopf map is given by 
$
  \pi_2 (z_1,z_2) = {z_1}/{z_2}. 
$
  \item {\it Orthogonal projection} $\pi_3 :{\mathbb H}^3 \to {\mathbb H}^2$  
defined as follows: Consider the Poincar\'e disc model of 
${\mathbb H}^3$ i.e. ${\mathbb H}^3$ identified with 
${\bf D}^3 = \left\{x= (x_1,x_2,x_3):\|x\| <1 \right\}$ with the 
metric ${4\sum d x_i^2}/{(1-|x\|^2)^2}$. Identify ${\mathbb H}^2$, 
via the Poincar\'e model, with the equatorial disc $x_3 =0$. Then 
$\pi_3$ is defined by projecting each $x\in {\mathbb H}^3$ to 
${\mathbb H}^2$ along the unique hyperbolic geodesic through $x$ which meets 
 ${\mathbb H}^2$ orthogonally. 
  \item {\it Orthogonal projection} 
$\pi_4 :{\mathbb H}^3 \to {\mathbb C}$ {\it to the plane at infinity} 
defined by considering the upper half-space model of 
${\mathbb H}^3 = \left\{ (x_1,x_2,x_3):x_3 >0 \right\}$ 
with the metric ${\sum d x_i^2}/{x_3}^2$ and 
setting $\pi_4 (x_1,x_2,x_3) = x_1 + i x_2$. 
\end{enumerate}
By Proposition~\ref{prop:ind conf}, the fibres of a non-constant 
\hm\ from a three-\dmen\ manifold $M^3$ form a \conf\ \cf\ by 
geodesics of $M^3$; the \conf s corresponding to the \hms\ 
defined in above examples are, respectively, 
\begin{enumerate}
  \item $\cf_1$, the foliation of ${\mathbb R}^3$ by vertical straight
lines;  
  \item $\cf_2$, the foliation of ${\mathbb S}^3$ by great circles
given by the intersection of $1$-dimensional complex subspaces of
${\mathbb C}^2$ with ${\mathbb S}^3 \subset {\mathbb C}^2$; 
  \item $\cf_3$, the foliation of ${\mathbb H}^3$ by geodesics
orthogonal 
to the equatorial disc, in the Poincar\'e disc model of ${\mathbb H}^3$;
  \item $\cf_4$, the foliation of ${\mathbb H}^3$ by vertical half 
lines, in the upper half-space model of ${\mathbb H}^3$. 
\end{enumerate}
Note that $\cf_1$, and $\cf_2$ are, in fact, Riemannian foliations. 
Now we can give the classification result for \hms\ from 
a 3-\dmen\ \scon\  space form \spf\ to a surface. 
\begin{thm}\cite{BaiWoo91}
\label{3dsc}
Up to isometries of \spf, a non-constant \hm\ of a 
three-\dmen\ \scon\  space form \spf\ to a surface $N^2$ is 
one of the examples $\pi_i$ described above followed by a 
weakly conformal map to $N^2$ and the associated \conf\ by 
geodesics of \spf\ is  one of the examples 
$\cf_i$ described above.  
\end{thm}

\subsection{The smoothing process of Baird and Wood} \qquad 

The smoothing process of Baird and the second author, obtained in \cite{BaiWoo92A}, 
is that an orbifold $O$, which is the leaf space of a Seifert 
fibre space without reflections can 
be smoothed to have a conformal structure.

Specifically, let $\left( M^3,\cf \right)$ be a Seifert fibre space 
 without reflections, with a $C^\infty$ 
metric such that the Seifert foliation 
\cf\ is a conformal foliation by geodesics. Let 
$L^M$ denote the leaf space of $M^3$; this is an orbifold with 
only cone points as singular points such that the cone points 
correspond to the critical fibres, cf.~\cite{scot83}. 
The horizontal conformality of \p\ at points on regular fibres 
gives $L^M \setminus \{{\rm cone\ points}\}$ a conformal 
structure $c$; then away from the critical fibres, the natural
projection \p\ of $\left( M^3,\cf \right)$
onto $(L^M \setminus \{{\rm cone\ points}\}, c)$ is a \hm. We have
\begin{prop}\cite[Prop.~2.6]{BaiWoo92A}
\label{prop:smooth}
With the above notation, there exists a unique $C^\infty$ 
conformal structure on $L^M$, such that if we denote 
$L^M$ with this conformal structure as $L^{M}_{s}$, the inclusion  
map $I:\left(
 L^M \setminus \{cone\ points\},c
\right) \to L^{M}_{s}$ is conformal. The composition 
$I\circ\p:M^3 \to L^{M}_{s}$ is a smooth \hm.
\end{prop}

We need a variant of this which we can state without mention of 
a 3-manifold:

A {\it 2-dimensional Riemannian orbifold\/} is an orbifold
which is locally 
the quotient of a 2-dimensional Riemannian manifold by a discrete group 
\Gp\ of isometries \cite{rat,scot83}, it has only cone points 
as singularities if and only 
if for each $p$, stabilizers $\Gp_p$ are all finite subgroups of 
orientation preserving isometries. 

Such an orbifold has a Riemannian metric and so a conformal 
structure away from the cone points. Then in the same way as 
\cite[Prop.~2.6]{BaiWoo92A} we have
\begin{prop}
Let $L$ be a 2-dimensional Riemannian orbifold with only cone 
points as singularities. Then there exists smooth conformal structure 
on $L$ such that, if we denote $L$ with this conformal 
structure by $L_s$, the identity map $I:L\to L_s$ is conformal 
on $L \setminus \{cone\ points\}$. 
\end{prop}

We shall call $L_s$ the {\it smoothed orbifold (associated to $L$)}.

%%
%%
%%
%%
%%
%%%%%%%%%%%%%%%%%%%%%%%%%%%%%%%%%%%%%%%%%%%%%%%%%%%%%%%%%%%
%%								
%%
%%		NEXT SECTION
%%
%%       H.Morph from 3-dim nonSC space forms
%%
%%%%%%%%%%%%%%%%%%%%%%%%%%%%%%%%%%%%%%%%%%%%%%%%%%%%%%%%%%%

\section{Harmonic morphisms from three-\dmen\ \nscon\ 
Euclidean and spherical space forms to a surface}
\label{sec:nscon}

\setcounter{equation}{0}

Throughout this section ${\mathbb E}^{m}_{1}$ will denote 
${\mathbb R}^m$ and ${\mathbb E}^{m}_{2}$ will denote 
${\mathbb S}^m$. 
Further  \G\ will denote a discrete group of isometries 
of ${\mathbb E}^{3}_{i}$  acting freely on ${\mathbb E}^{3}_{i}$. 
Note that such a \G\ acts properly discontinuously \cite[p~406]{scot83} and 
the quotient $M^3={\mathbb E}^{3}_{i} /\G$ is a Euclidean or 
spherical space form, i.e. a connected complete 3-dimensional 
\rman\ of constant non-negative curvature \cite[p~69]{wol66}. 
Conversely any Euclidean or spherical 3-dimensional space form 
is homothetic to such a quotient. 

We first describe some {\it standard} \hms\ from such a space form 
to a surface:
\begin{thm}
\label{thm:std}
Let $\pi_i :{\mathbb E}^{3}_{i}\to {\mathbb E}^{2}_{i}$ $(i=1,2)$ 
be one of the standard \hms\ in \S\ref{sec:531} and 
$\cf_i$ be the corresponding Riemannian foliation by geodesics. 
Suppose that \G\ is a discrete group of isometries acting 
freely on ${\mathbb E}^{3}_{i}$ such that 
\begin{description}
  \item[(a)]\label{cond a} \G\ preserves $\cf_i$, i.e. \G\ 
   maps leaves of $\cf_i$ to leaves.\\
Then 
 \begin{description}
        \item[(i)] \G\ descends through $\pi_i$ to an action of 
         ${\mathbb E}^{2}_{i}$ by a group \Gp\ of isometries, 
        \item[(ii)] \cf\ factors to a Riemannian foliation 
         $\cf_{i,\G}$ by geodesics 
         of the space form $M^3={\mathbb E}^{3}_{i} /\G$ with leaf space 
         $L^M={\mathbb E}^{2}_{i} /\Gp$. 
 \end{description}
Suppose further that 
  \item[(b1)]\label{cond b1} for any $p\in {\mathbb E}^{2}_{i}$ 
  the stabilizer $\Gp_p \subset \Gp$ of $p$ contains no reflections,
  \item[(b2)]\label{cond b2} \Gp\ acts discontinuously on 
    ${\mathbb E}^{2}_{i}$. \\
Then  
\begin{description}
        \item[(iii)] $L^M={\mathbb E}^{2}_{i} /\Gp$ is a Riemannian 
         orbifold whose only possible singularities are cone 
         points;
        \item[(iv)] Letting  $L^{M}_{s}$ denote the smoothed orbifold and 
    $I:L^M\to L^{M}_{s}$ the identity map, $\pi_i$ factors 
    to a continuous map $\pi_{i,\G}:M^3\to L^M$ such that the composition 
 $\p_{i,\G}:M^3\mathop {\longrightarrow} \limits^{\pi_{i,\Gamma}} L^M 
 \mathop {\longrightarrow} \limits^I L^{M}_{s}$ 
is a smooth \hm. We thus have a 
    commutative diagram, where the vertical arrows are natural projections. 
 \end{description}
\end{description}
\begin{figure}[htbp]
\begin{picture}(210,150)(-120,0)
\put(40,130){\makebox(0,0){${\mathbb E}^{3}_{i}$}}
\put(50,130){\vector(1,0){70}}
\put(85,140){\makebox(0,0){$\hat{\pi_{i}}$}}
\put(130,130){\makebox(0,0){${\mathbb E}^{2}_{i}$}}
\put(130,120){\vector(0,-1){60}}
%\put(50,90){\makebox(0,0){$\pi$}}
\put(30,50){\makebox(0,0){$M^{3}={\mathbb E}^{3}_{i}/\Gamma$}}
\put(40,120){\vector(0,-1){60}}
%\put(120,90){\makebox(0,0){$\pi ^\prime$}}
\put(140,50){\makebox(0,0){$L^{M}={{\mathbb E}^{2}_{i}}/\Gamma ^{\prime}$}}
\put(60,50){\vector(1,0){45}}
\put(85,60){\makebox(0,0){$\pi_{i,\G}$}}
\put(195,60){\makebox(0,0){$I$}}
\put(170,50){\vector(1,0){50}}
\put(235,50){\makebox(0,0){$L^{M}_{s}$}}
\put(130,12){\makebox(0,0){$\phi_{i,\G}$}}
\qbezier(40,43)(115,5)(220,43)
\put(220,43){\vector(3,1){0}}
%
%\put(140,15){\makebox(0,0){$L^{M}_{s}$}}
%\put(140,40){\vector(0,-1){15}}
%\put(130,32){\makebox(0,0){$id$}}
%\put(191,-27){\makebox(0,0){$L^{M}_{s}$}}
%\put(140,120){\vector(1,-2){75}} 
%\put(180,60){\makebox(0,0){$\hat{\xi}$}}
%\put(50,40){\vector(2,-1){125}}
%\put(130,-10){\makebox(0,0){$\phi_{i,\G}$}}
%\put(150,40){\vector(1,-2){30}}
%\put(174,7){\makebox(0,0){$I$}}
%\put(30,130){\line(-1,0){90}}
%\put(-60,130){\line(0,-1){157}}
%\put(-60,-27){\vector(1,0){237}}
%\put(-70,52){\makebox(0,0){$\hat{\phi}$}}
%\put(150,10){\vector(1,-1){60}}
%\put(175,0){\makebox(0,0){$\xi$}}
\end{picture}
%\caption{ }\label{dia:con compt}
\end{figure}
\end{thm}
\begin{rk}\label{rem:std}
\hfill
\begin{enumerate}
\item We shall see later that Conditions~{\bf (b1)},~{\bf (b2)} 
are {\it necessary\/} for 
{\bf (iii)} and {\bf (iv)}. We shall again call the \hms\ $\phi_{i,\G}$ 
{\it standard \hms}. 
\item The same proposition holds with minor changes for the two 
standard \hms\ from ${\mathbb H}^3$. However, as quotients of ${\mathbb H}^3$
have not yet been classified, we shall not consider this case further.
(Note also that there are no non-constant harmonic morphisms from
compact quotients of ${\mathbb H}^3$ \cite{Bai90}.)
\end{enumerate}
\end{rk}
\begin{proof}
\hfill 
\begin{description}
        \item[(i)] Since \G\ preserves $\cf_i$ it induces an action on
  the leaf space ${\mathbb E}^{2}_{i}$; since 
  \G\ preserves the distance between leaves of the Riemannian foliations 
  $\cf_i$, this action is by a group $\Gp$ of isometries.
        \item[(ii)] This is clear since the natural projection 
  ${\mathbb E}^{3}\to {\mathbb E}^{3}/\G$ is a Riemannian covering. 
        \item[(iii)] Immediate from {\bf (i)} by definition of Riemannian 
  orbifold.
        \item[(iv)] By definition of \Gp, $\pi_i$ factors to a continuous map 
  $\pi_{i,\G}:M^3\to L^M$ which is an isometry away from cone points. 
  Since the identity map $I:L^M\to L^{M}_{s}$ is conformal away from cone 
points, the composition $\p_{i,\G}=I\circ \pi_{i,\G}:M^3\to L^{M}_{s}$ is 
$C^0$ and is a $C^\infty$ \hm\ on 
$M^3\setminus \pi^{-1}_{i,\Gamma}({\rm cone points})$. 
As in \cite[Theorem~2.18]{BaiWoo91} it is a $C^\infty$ \hm. 
Since the composition  
${\mathbb E}^{2}\to L^M ={\mathbb E}^{2}/\Gp\mathop \to \limits^I L^{M}_{s}$ 
has branch points  
precisely at the inverse images of cone points, $\p_{i,\G}$ has critical 
points precisely at $\pi^{-1}_{i,\G}({\rm cone points})$. 
\end{description}
\end{proof}
{\bf Remarks.}\\
If $M^3$ is compact, $(M^3,\cf_{i,\G})$ is a Seifert fibre space without 
reflections and our construction is the same as that in 
\cite[Prop.~2.6]{BaiWoo92A}. 
However, if $M^3$ is not compact, the leaves of $\cf_{i,\G}$ may 
not be closed and our construction is a little more general.\\
%\medskip
{\bf Example~1.}\\
Suppose that \G\ is the infinite cyclic group generated by the glide 
reflection 
\linebreak
$(x_1,x_2,x_3)\mapsto (-x_1,x_2,x_3 +1)$. This preserves the 
standard foliation $\cf_1$ on ${\mathbb R}^3$ and factors to the 
group \Gp, of order 2, of isometries of ${\mathbb R}^2$ generated by 
the reflection $(x_1,x_2)\mapsto (-x_1,x_2)$ so that 
${\mathbb R}^2 /\Gp$ is the half plane, an orbifold with reflector line 
$x_2 =0$. $M^3={\mathbb R}^3 /\G$ is the product of ${\mathbb R}$ and 
an infinite width Moebius strip and the foliation $\cf_1$ factors to 
a Seifert fibration $\cf_{1,\G}$ on $M^3$ with reflections, with critical 
fibres the vertical lines in the plane $x_2=0$. However, \G\ does not 
satisfy Condition~{\bf (b1)} and there is no \hms\ from 
$M^3$ with associated foliation $\cf_{1,\G}$.\\
{\bf Example~2.}\\
Suppose that \G\ is the infinite cyclic group generated by the 
screw motion $({\mathcal R}_\theta,t_{v_3})$ consisting of 
a rotation about the $x_3$-axis through an angle $\theta$
and a translation by one unit in the 
$x_3$-direction. Then \G\ preserves the standard foliation $\cf_1$ 
on ${\mathbb R}^3$ and factors to the cyclic 
group \Gp\ of isometries of ${\mathbb R}^2$ generated by rotation 
${\mathcal R}_\theta$ through $\theta$ about the origin. There are 
two cases:
\begin{enumerate}
\item If $\theta/2\pi$ is irrational, \Gp\ is infinite cyclic, 
${\mathbb R}^3 /\G$ is diffeomorphic to ${\mathbb R}^3$ and 
$\cf_1$ factors to a foilation $\cf_{1,\G}$ of 
${\mathbb R}^3 /\G$ by straight lines. \G\ does not satisfy 
Condition~{\bf (b2)} and there is no \hm\ from ${\mathbb R}^3 /\G$ 
to a surface with corresponding  foliation $\cf_{1,\G}$.
\item If  $\theta/2\pi$ is rational, say $p/q$ in lowest terms, 
\Gp\ is cyclic of order $q$ generated by rotation through $2\pi/q$, 
$L^M={\mathbb R}^2 /\Gp$ is an orbifold with one cone point of angle 
$2\pi/q$, \ ${\mathbb R}^3 /\G$ is diffeomorphic to 
${\mathbb R}^2 \times {\mathbb S}^1$ and $\cf_1$  factors to a 
Seifert fibration $\cf_{1,\G}$ of ${\mathbb R}^3 /\G$ by 
circles. $L^{M}_{s}$ can be identified with ${\mathbb C}$ via the 
homeomorphism ${\mathbb R}^2 /\Gp\to {\mathbb R}^2 ={\mathbb C}$,  
$[z]\mapsto z^q$ and the resulting standard \hm\ is given by 
$$
{\mathbb R}^3 /\G \to {\mathbb R}^2, \ \ [(z,t)]\mapsto z^q .
$$
\end{enumerate}
%\qquad \\
\medskip
We next show that, up to postcomposition with weakly conformal maps, 
our standard \hms\ give {\it all\/} \hms\ from complete flat or spherical 
3-dimensional space forms. Let $\GG_i$ denote the set of all discrete 
groups of isometries acting freely on ${\mathbb E}^{3}_{i}$ and 
satisfying {\bf (a)}, {\bf (b1)} and {\bf (b2)} of Theorem~\ref{thm:std}.
\begin{thm}
Let $M^3$ be a complete 3-dimensional space form of non-negative 
curvature and let $\p:M^3\to N^2$ be a non-constant \hm. Then 
\begin{enumerate}
\item $M^3$ is homothetic to ${\mathbb E}^{3}_{i} /\G$ for 
some $\G\in \GG_i$
\item \p\ is the composition of a homothety 
$M^3\to {\mathbb E}^{3}_{i} /\G$, a standard \hm\ 
$\p_{i,\G}:{\mathbb E}^{3}_{i} /\G \to L^{M}_{s}$ and a weakly 
conformal map $L^{M}_{s}\to N^2$.
\end{enumerate}
\end{thm}
\begin{proof}
As in \cite[p~69]{wol66} $M^3$ is homothetic to 
${\mathbb E}^{3}_{i} /\G$ for some discrete group \G\ of 
isometries acting freely on ${\mathbb E}^{3}_{i}$. Identifying 
$M^3$ with ${\mathbb E}^{3}_{i} /\G$, the composition 
$\hat{\p}:{\mathbb E}^{3}_{i}\to {\mathbb E}^{3}_{i} /\G\to N^2$ 
is a \hm. By \cite{BaiWoo88} (see Theorem~\ref{3dsc} above) after 
applying an isometry of ${\mathbb E}^{3}_{i}$, this is the composition 
of the standard \hm\ $\pi_i :{\mathbb E}^{3}_{i} \to {\mathbb E}^{2}_{i}$ 
with a weakly conformal map $\hat{\xi}:{\mathbb E}^{2}_{i} \to N^2$. 
\G\ preserves the corresponding foliation $\cf_i$ which therefore 
descends to a Riemannian foliation $\cf_{i,\G}$ by geodesics of 
${\mathbb E}^{3}_{i} /\G \simeq M^3$. Thus \G\ satisfies 
Condition~{\bf (a)}.

Further \Gp\ acts by isometries on ${\mathbb E}^{2}_{i}$ and, 
by diagram, $\hat{\xi}$ factors to a continuous map 
$\xi:{\mathbb E}^{2}_{i} /\Gp \to N^2$ conformal away from critical 
points so that we have a commutative diagram.  

\begin{figure}[htbp]
\begin{picture}(210,185)(-100,-10)

\put(40,130){\makebox(0,0){${\mathbb E}^{3}_{i}$}}
\put(50,130){\vector(1,0){70}}
\put(85,140){\makebox(0,0){$\hat{\pi _{i}}$}}

\qbezier(50,140)(140,180)(236,96)
\put(240,94){\vector(2,-1){0}}

\put(130,130){\makebox(0,0){${\mathbb E}^{2}_{i}$}}
\put(130,120){\vector(0,-1){60}}
\put(50,90){\makebox(0,0){$\pi$}}

\put(30,50){\makebox(0,0){$M^3={\mathbb E}^{3}_{i}/\G$}}
\put(40,120){\vector(0,-1){60}}
\put(120,90){\makebox(0,0){$\pi ^\prime$}}

\qbezier(50,40)(140,0)(237,78)
\put(240,80){\vector(2,1){0}}

\put(140,50){\makebox(0,0){$L^{M}={\mathbb E}^{2}_{i}/\Gamma ^{\prime}$}}
\put(60,50){\vector(1,0){45}}
\put(85,60){\makebox(0,0){$\pi _{i,\G}$}}

%\put(140,15){\makebox(0,0){$L^{M}_{s}$}}
%\put(140,40){\vector(0,-1){15}}
%\put(130,32){\makebox(0,0){$id$}}

\put(250,90){\makebox(0,0){$N^2$}}
\put(145,130){\vector(2,-1){83}} 
\put(175,105){\makebox(0,0){$\hat{\xi}$}}

%\put(50,40){\vector(2,-1){155}}
\put(130,20){\makebox(0,0){$\phi$}}

\put(130,165){\makebox(0,0){$\hat{\phi}$}}

\put(172,55){\vector(2,1){59}}
\put(185,70){\makebox(0,0){$\xi$}}

\end{picture}
%\caption{A Commutative Diagram}\label{F_is_smooth}
\end{figure}

Now, if $\Gp$ did not act discontinuously, then, in the neighbourhood of
some $p \in M^3$ there would be an infinite sequence of line segments
mapped by $\phi$ to the same point, which is impossible by the local
form of a \hm\ (see \S \ref{prop:ind conf}), thus Condition {\bf (b2)}
is satisfied..

Further for any $p\in {\mathbb E}^{2}_{i}$ the action of its stabilizer 
$\Gp_p$ preserves the fibres of $\hat{\xi}$, since any point is either a 
regular point of $\hat{\xi}$ or a branch point of order $k\geq 2$, this means 
$\Gp_p$ lies in a cyclic subgroup of order $k\geq 1$ of $SO(2)$ 
so that Condition {\bf (b1)} is satisfied.  We have shown that $\G \in
\GG_1$.

Then ${\mathbb E}^{2}_{i} /\Gp$ is 
an orbifold with only cone points. Note that it is the leaf space 
$L^M$ of $\cf_{i,\G}$, it can thus be smoothed as above and the 
 resulting \hm\ 
$\p_{i,\G}:{\mathbb E}^{3}_{i} /\G \to L^M\mathop \to \limits^I L^{M}_{s}$ is 
a standard \hm. Clearly $\xi:L^{M}_{s}\to N^2$ is $C^0$, and 
$C^\infty$  and conformal away from the cone points, therefore 
conformal everywhere and the Theorem is proved.
\end{proof}
%%%%%%%%%%%%%%%%%%%%%%%%%%%%%%%%%%%%%%%%%%%%%%%%%%%%%%%%%%%%%%%%
%%
%%	next section
%%
%%
%%%%%%%%%%%%%%%%%%%%%%%%%%%%%%%%%%%%%%%%%%%%%%%%%%%%%%%%%%%%%%%%
\section{the standard harmonic morphisms from complete flat and 
spherical space forms}\label{sec:standard}
To complete our description of \hms\ from Euclidean and spherical space 
forms $M^3$ we list the standard \hms\ for these two cases. 

%First the Euclidean case.
\subsection{The Euclidean case}
Recall that $\GG_1$ denotes the set of all discrete 
groups of isometries acting freely on ${\mathbb R}^{3}$ and 
satisfying {\bf (a)}, {\bf (b1)} and {\bf (b2)} of Theorem~\ref{thm:std}. 
\begin{thm}
%\hfill
Let $M^3$ be a complete flat space form.
\begin{enumerate}
\item Any \hm\ from $M^3$ to a surface $N^2$ is the composition
of an isometry  $M^3\to {\mathbb R}^3 /\G$, a standard \hm\
$\p_{1,\Gamma}$ for some $\G \in \GG_1$ and a
weakly conformal map $L^{M}_{s} \to N^2$.  
\item The standard \hms\ are given by 
$$
\p_{1,\Gamma}=I\circ \pi_{1,\Gamma}:{\mathbb R}^3 /\G
 {\longrightarrow} {\mathbb R}^2 /\Gp 
 =L^M\mathop {\longrightarrow} \limits^I L^{M}_{s}
$$
where $\Gamma \in \GG_1$.

\item For each isometry type of $M^3$, there are (non-constant) standard
harmonic morphisms.  We list these for orientable $M^3$.  In the
following list, with notations explained below,
for each of the 10 diffeomorphism types A of $M^3$, we
list the possible groups $\G$ and the parameters they depend on,
then the different choices (a), (b), ... of parameters which ensure that
$\G \in \GG_1$.  For each such choice, we list the corresponding $\Gp$
and $L^M$.
\end{enumerate}

$M^3$ non-compact:

\begin{enumerate}
\item $A = {\mathcal E}$, $M^3 \approx {\mathbb R}^3$, $\Gamma = \{e\}$:
	\begin{enumerate}
        \item $\Gp=\{e\}$, $L^M={\mathbb R}^2$.
        \end{enumerate}
\item $A = {\mathcal J}^{\theta}_{1}$, $M^3
	\approx {\mathbb S}^1 \times {\mathbb R}^2$,
	$\Gamma=\la (R_{\theta}(v), t_v \ra)$ where
	$\theta \in {\mathbb R}$ and $v \in {\mathbb R^3}$ is non-zero:
	\begin{enumerate}
	\item $v \parallel e_3$, $\theta = 0$: $\Gp = \{e\}$, $L^M={\mathbb R}^2$.
	\item $v \parallel e_3$, $\theta = \pi$: $\Gp = \la R_{\pi} \ra$,
      		$L^M= {\mathbb R}^2(2) =$ cone of angle $\pi$.
	\item $v \parallel e_3$, $\theta\ne 0,\pi$ and $\theta=2\pi p/q$ where 
       		$p,q\in {\mathbb Z}, \ (p,q)=1, \ q\ne0$: $\Gp = \la
		R_{\theta} \ra$, $L^M={\mathbb R}^2 (q)$ = cone of angle $2\pi/q$.
	\item $v \not\parallel e_3$, $\theta = 0$: $\Gp =
		\la t_{\pi(v)} \ra$, $L^M =$cylinder.
	\item $v \in e_{3}^{\perp}$, $\theta = \pi$: $\Gp = \la
		(S_v, t_v) \ra$,
      		$L^M=$M\"obius band.
	\end{enumerate}
\item $A = {\mathcal T}_1$, $M^3 \approx
	{\mathbb S}^1 \times {\mathbb S}^1 \times {\mathbb R}$,
	 $\Gamma =\la t_{v_1}, t_{v_2} \ra$ where $v_1, v_2 \in {\mathbb
R}^3$ are linearly independent:
	\begin{enumerate}
	\item $e_3 \notin \spn \{v_1,v_2\}$: $\Gp = \la t_{\pi(v_1)},
		t_{\pi(v_2)} \ra$, $L^M=$ 2-torus.
	\item $e_3 \in \spn \{v_1,v_2\}$ and $\{ \pi(v_1)$, $\pi(v_2) \}$
		rationally related:
    		$\Gp = \la t_w \ra$ where 
		$\spn_{\mathbb Z}\{w\} = \spn_{\mathbb Z}\{\pi(v_1)$, $\pi(v_2)\}
		\in {\mathbb R}^2$, $L^M=$ cylinder. 
	\end{enumerate}
\item  $A = {\mathcal K}_1$, 
	$\G = \la t_{v_1}, (R_{\pi}(v_2),t_{v_2}) \ra$
	where $v_1, v_2 \in {\mathbb R}^3$ are linearly independent. 
	\begin{enumerate}
	\item $v_1\not\parallel e_3$, $v_2\in e_{3}^{\perp}$:
    		$\Gp: = \la t_{\pi(v_1)}, (S_{v_2},t_{v_2})
		\ra$, $L^M=$Klein bottle.
	\item $v_1\parallel e_3$, $v_2 \in e_3^{\perp}$,
	$\Gp = \la (S_{v_2}, t_{v_2}) \ra$, $L^M=$M\"obius band.
	\item $v_2\parallel e_3$:
    	$\Gp = \la t_{\pi(v_1)}, R_{\pi} \ra$, $L^M= D^2(2,2)$.
	\end{enumerate}
\bigskip
$M^3$ compact:

\item $A = {\mathcal G}_1$, $M^3 =$ 3-torus,
	$\G = \la t_{v_1},t_{v_2},t_{v_3} \ra$ where $v_1, v_2, v_3 \in
{\mathbb R}^3$ are linearly independent:
        \begin{enumerate}
	\item $\{ \pi(v_1), \pi(v_2), \pi(v_3) \}$ rationally related: $\Gp = \la
		t_{w_1}, t_{w_2} \ra$ where $\spn_{\mathbb Z}\{w_1,w_2 \} =
		\spn_{\mathbb Z}\{ \pi(v_1), \pi(v_2), \pi(v_3) \}$,
		$L^M$=2-torus.
	\end{enumerate}
\item $A = {\mathcal G}_2$, $\G =  \la (R_{\pi}(v_1),t_{v_1/2}),
	t_{v_1}, t_{v_2}, t_{v_3} \ra$ where $v_1,v_2,v_3 \in
{\mathbb R}^3$ are linearly independent and $v_1 \perp \spn \{v_2,v_3\}$: 
	\begin{enumerate}
	\item $v_1\in e_{3}^{\perp}$, $\{\pi(v_2), \pi(v_3) \}$ rationally related:
  		$\Gp = \la (S_{v_1},v_1), t_w \ra$ where $\spn_{\mathbb Z}\{w\} =
		\spn_{\mathbb Z}\{\pi(v_2), \pi(v_3)\}$, $L^M$=Klein bottle.
	\item $v_1\parallel e_{3}$: $\Gp =
		\la R_{\pi}, t_{v_2}, t_{v_3} \ra$,
		$L^M={\mathbb S}^{2}(2,2,2,2)$.
	\end{enumerate}
\item $A = {\mathcal G}_3$, $\G = \la (R_{2\pi/3}(v_1), t_{v_1/3}),
	t_{v_2}, t_{v_3}, t_{v_3} \ra$  where  $v_1,v_2,v_3 \in
	{\mathbb R}^3$ are linearly independent, $v_1 \perp
	\spn\{v_2,v_3\}$, $\|v_2\|=\|v_3\|$, $\angle(v_2,v_3) = 2\pi/3$:
	\begin{enumerate}
	\item $v_1\parallel e_{3}$:
  		$\Gp = \la R_{2\pi/3}, t_{v_2}, t_{v_3} \ra$,
                $L^M={\mathbb S}^{2}(3,3,3)$.
	\end{enumerate}
\item $A ={\mathcal G}_4$, $\G = \la (R_{\pi/2}(v_1), t_{v_1/4}),
	t_{v_1}, t_{v_2}, t_{v_3} \ra$ where $v_1,v_2,v_3 \in {\mathbb
	R}^3$ are mutually orthogonal with $\|v_2\| = \|v_3\|$:
	\begin{enumerate}
	\item $v_1\parallel e_{3}$:
  		$\Gp = \la R_{\pi/2}, t_{v_2}, t_{v_3} \ra$,
  		$L^M={\mathbb S}^{2}(2,4,4)$.
	\end{enumerate}
\item $A = {\mathcal G}_{5}$, $\G = \la (R_{\pi/6}(v_1), t_{v_1/6}),
	t_{v_1}, t_{v_2}, t_{v_3} \ra$ where $v_1,v_2,v_3 \in {\mathbb
	R}^3$ are linearly independent, $v_1 \perp \spn\{ v_2, v_3 \}$,
	$\|v_2\|=\|v_3\|$ and $\angle(v_2,v_3) = \pi/6$:	
	\begin{enumerate}	
	\item $v_1\parallel e_{3}$:
  		$\Gp = \la R_{\pi/6},  t_{v_2}, t_{v_3} \ra$,
  		$L^M={\mathbb S}^{2}(2,3,6)$.
  	\end{enumerate}	
\item $A = {\mathcal G}_{6}$, $\G = \la(R_{\pi}(v_1), t_{v_1/2}),
	(R_{\pi}(v_2), t_{(v_2+v_3)/2}), (R_{\pi}(v_3),
	t_{(v_1+v_2+v_3)/2}), t_{v_1}, t_{v_2}, t_{v_3} \ra$ where
	$v_1,v_2,v_3$ are mutually orthogonal:
 	\begin{enumerate}
	\item $v_1 \parallel e_{3}$:
  		$\Gp = \la R_{\pi},
		(S_{v_2}, t_{(v_2+v_3)/2}), (S_{v_3},t_{(v_2+v_3)/2)}),
		 t_{v_2}, t_{v_3} \ra$, $L^M = {\mathbb P}^2 (2,2)$.
	\item $v_2\parallel e_{3}$:
		$\Gp = \la (S_{v_1}, t_{v_1/2}), (R_{\pi},t_{v_3/2}),
		(S_{v_3}, t_{(v_1+v_3)/2}), t_{v_1}, t_{v_3} \ra$,
		$L^M = {\mathbb P}^2 (2,2)$.
 	\item $v_3\parallel e_{3}$:
        	$\Gp = \la (S_{v_1}, t_{v_1/2}), (S_{v_2}, t_{v_2/2}),
		(R_{\pi}, t_{(v_1+v_2)/2}), t_{v_1}, t_{v_2} \ra$,
		$L^M = {\mathbb P}^2 (2,2)$.
	\end{enumerate}
\end{enumerate}
\end{thm}

\noindent {\it Notation} \\
For $v \in {\mathbb R}^3$, $t_v$ denotes
translation in ${\mathbb R}^3$ through $v$, and for $v \neq 0$, $\theta
\in \mathbb R$, $R_{\theta}(v)$
denotes rotation in ${\mathbb R}^3$ through an angle $\theta$ about $v$.

For $v \in {\mathbb R}^2$, $t_v$ denotes translation in ${\mathbb
R}^2 = e_3^{\perp}$ through $v$,
for $v \neq 0$, $S_v$ denotes reflection in $v$ and, for $\theta \in
{\mathbb R}$, $R_{\theta}$ denotes rotation
about the origin through an angle $\theta$.
 
Bracketted pairs $(A,t)$ denote $A$ followed by $t$, thus, for example,
$(R_{\theta}(v), t_v)$ is the screw motion through $\theta$ along $v$.

We say that vectors $v_i$, $(i=1,\cdots, k)$ are {\em rationally related\/}
if they are linearly dependent over $\mathbb Q$.  If $k =2$ and the
$v_i$ are not all parallel, this means that the corresponding
translations $t_{v_i}$ generate a $2$-dimensional lattice.

$N(r,s, \cdots)$ denotes the orbifold with underlying surface $N$ and
cone points of orders $r,s,\cdots$.
  
\begin{proof}
\hfill

Parts 1. and 2. follow from Theorems 3.1 and 3.3 noting that we can
scale the homothety of $M^3$ to ${\mathbb R}^3$ to make it an isometry.
For Part 3, the discrete groups  \G\ acting freely on ${\mathbb R}^3$
and giving orientable quotients ${\mathbb R}^3/\G$ are classified 
in \cite[Theorem 3.5.1, 3.5.5]{wol66} up to affine diffeomorphism, we add
parameters according to the remarks at the end of \cite[\S 3.5]{wol66} to
obtain all groups \G\ up to translational equivalence.
Then we determine
the values of those parameters for which $\G\in \GG_1$,
using the following simple facts:  (i) condition (b1) is automatic for ${\mathbb
R}^3/\G$ orientable, (ii) a rotation through an angle $\theta \in (0,2\pi)$ about
a vector $v$ preserves the $e_3$-direction if and only if either $v \parallel e_3$ or
$v \perp e_3$ and $\theta = \pi$. The induced group $\Gp$ and the
resulting orbifold $L^M = {\mathbb R}^2 / \Gp$ are then calculated.
\end{proof}

\noindent {\bf Remarks.}

1. Each of the 10 diffeomorphism types represents exactly one {\it
affine\/} diffeomorphism type of $M^3$ except for ${\mathcal J}^{\theta}_{1}$,
where the affine diffeomorphism type is parameterized by $\theta \mod \pi$.

2. For 1, 2(d), 2(e), 3(a), 4(a) the foliation $\cf_{1,\G}$ defined by
$\G$ has fibres which are lines, for the other cases they are circles
and $(M^3, \cf_{1,\G})$ is a Seifert fibre space.
  
3. A similar list can be given for $M^3$ non-orientable.

\subsection{The spherical case} \quad
\newline
First note that, since orientation reversing 
transformations have fixed points, any discrete subgroup of 
${\rm Isom}({\mathbb S}^3)=O(4)$ acting freely is a finite 
subgroup of $SO(4)$.

Identify ${\mathbb R}^4$ with the quaternions and ${\mathbb S}^3$ 
with the unit quaternions. Let $\psi:{\mathbb S}^3\to SO(3)$ be the standard
double cover given by $\psi(q) = a\mapsto qaq^{-1}$ and let 
$\p:{\mathbb S}^3 \times {\mathbb S}^3\to SO(4)$ be the double cover 
given by $\p(q_1,q_2) =x\mapsto q_1 xq_2^{-1}$. Let
$\G_1 = \p({\mathbb S}^1 \times {\mathbb S}^3)$ and
$\G_2 = \p({\mathbb S}^3 \times {\mathbb S}^1)$.
Let $p:SO(4)\to SO(3)\times SO(3)$ be the unique map 
defined by 
$\psi \times \psi = p\circ \p$ so that $p$ is a surjective homomorphism with
Ker~$p=\{I,-I\}\subset SO(4)$, 
and let $H_i\subset SO(3)$ $(i=1,2)$ be the projections of 
$H=p(\G)$ onto the two factors.
Then a finite subgroup $\G$
of $SO(4)$ preserves  the standard Hopf foliation $\cf_2$ if and only if
either $\G \subset \G_1$, or $\G \subset \G_2$ with $H_1$ dihedral with
cyclic subgroup in $\psi(S^1)$ and
$H_2$ is cyclic  \cite{scot83}.  In the first case $H_1$ is cyclic and 
$H_2$ is cyclic, dihedral or the symmetry group of a regular solid
namely, the tetrahedral, octahedral or icosahedral groups $\bf T$, $\bf
O$ or $\bf I$. Further 
the induced action of $(q_1,q_2) \in {\mathbb S}^1 \times {\mathbb S}^3$
on the leaf space ${\mathbb S}^2$ 
is given by $\psi(q_2)$ so that $L^M = S^2/H_2$.
In the second case, we must additionally factor out by the action of an
element of order $2$ in $H_1$ giving an orbifold with underlying space
real projective $2$-space ${\mathbb P}^2$.  Thus we have 
\begin{thm}
%\hfill
Let $M^3$ be a spherical space form. 
\begin{enumerate}
\item Any \hm\ from a complete spherical space form $M^3$ to a surface 
$N^2$ is the composition of a homothety $M^3\to {\mathbb S}^3 /\G$, 
a standard \hm\ $\p_{2,\G}$ for some finite subgroup $\G$ of $\G_1$
and a weakly conformal map  $L_{s}^{M}\to N^2$.
\item The standard  harmonic morphisms are given by
$$
\p_{2,\G}=I\circ\pi_{2,\G}:{\mathbb S}^3 /\G \longrightarrow
{\mathbb S}^2 /\Gp=L^M\mathop {\longrightarrow} \limits^I L^{M}_{s}
$$
where $\G$ is a finite subgroup of $\G_1$.
\item  For each isometry type of $M^3$ there are (non-constant) standard
harmonic morphisms.  In the following table, we categorize these
according to the subgroups $H_1$ and $H_2$; we list the corresponding quotient
$S^3/\G$ and orbifold $L^M$.
\end{enumerate}

%\pagebreak 
\begin{center}
%\begin{table}
\begin{tabular}{|l|l|l|l|}   \hline 
\ $H_1$\ &\ $H_2$\     &\ $M^3={\mathbb S}^3 /\G$ &\ $L^M$\ \\
\hline  
\ ${\mathbb Z}_p$\ &\ ${\mathbb Z}_q$\ &\ Lens spaces  &\ ${\mathbb S}^{2}(q,q)$  \\
\hline 
\ ${\mathbb Z}_p$\ &\ ${\bf D_{m}}$  &\ Prism spaces &\ ${\mathbb S}^{2}(2,2,m)$  \\
\hline
  ${\bf D}_m$\     &\ ${\mathbb Z}_q$  &\ Prism spaces &\ ${\mathbb P}^2(q)$ \\
\hline
\ ${\mathbb Z}_p$\ &\ ${\bf T}$  &\ Tetrahedral spaces\ &\ ${\mathbb S}^{2}(2,3,3)$ \\
\hline
\ ${\mathbb Z}_p$\ &\ ${\bf O}$  &\ Octahedral spaces  &\ ${\mathbb S}^{2}(2,3,4)$ \\
\hline
\ ${\mathbb Z}_p$\ &\ ${\bf I}$  &\ Icosahedral spaces\  &\ ${\mathbb S}^{2}(2,3,5)$ \\
\hline 
\end{tabular}
%\end{table}
\end{center}
\end{thm}

\begin{proof}

Parts 1. and 2. follow from Theorems 3.1 and 3.3.
 
For Part 3., that there are standard harmonic morphisms for any isometry type follows
from the fact that any finite subgroup of $SO(4)$ acting freely on 
${\mathbb S}^3$ 
is conjugate in $O(4)$ to a
subgroup of $\Gamma_1$ \cite{scot83}.
\end{proof}

\subsection*{Acknowledgments}
The first named author (MTM) would like to express his 
gratitude to 
The Director~I.C.T.P. for the hospitality and support.

%\bibliographystyle{amsplain}
%\bibliography{tahir}

\begin{thebibliography}{10}

\bibitem{Bai83}
P.~Baird, \emph{{H}armonic maps with symmetry, harmonic morphisms, and
  deformation of metrics}, Pitman Res.\ Notes Math.\ Ser., vol.~87, Pitman,
  Boston, London, Melbourne, 1983.

\bibitem{Bai90}
P.~Baird, \emph{ {H}armonic morphisms and circle actions on {$3$}- and
   {$4$}-manifolds}, Ann.\ Inst.\ Fourier (Grenoble) \textbf{40} (1990),   
    177--212.

\bibitem{BaiEel81}
P.~Baird and J. Eells, \emph{ {A} conservation law for harmonic maps}, 
{G}eometry {S}ymposium. {U}trecht 1980, 
Lecture Notes in Math., \textbf{894} (1981), 1--25.

\bibitem{BaiWoo88}
P.~Baird and J.C. Wood, \emph{{B}ernstein theorems for harmonic morphisms from
  {${\mathbb {R}}^3$} and {${\mathbb S}^3$}}, Math.\ Ann. \textbf{280} (1988), 579--603.

\bibitem{BaiWoo91}
\bysame, \emph{{H}armonic morphisms and conformal foliations by geodesics of
  three-dimensional space forms}, J. Austral.\ Math.\ Soc.\ Ser.\ A \textbf{51}
  (1991), 118--153.

\bibitem{BaiWoo92A}
\bysame, \emph{{H}armonic morphisms, {S}eifert fibre spaces and conformal
  foliations}, Proc.\ London Math.\ Soc. \textbf{64} (1992), 170--196.

\bibitem{Fug78}
B.~Fuglede, \emph{{H}armonic morphisms between {R}iemannian manifolds}, Ann.\
  Inst.\ Fourier (Grenoble) \textbf{28} (1978), 107--144.

\bibitem{GUD}
S.~Gudmundsson, \emph{On the existence of harmonic morphisms from symmetric 
      spaces of rank one}, preprint, Lund University (1996).

%\bibitem{GudA}
%S.~Gudmundsson, \emph{{H}armonic morphisms from complex projective spaces},
%  Geom.\ Dedicata \textbf{53} (1994), 155--161.
%
\bibitem{Ish79}
T.~Ishihara, \emph{{A} mapping of {R}iemannian manifolds which preserves
  harmonic functions}, J. Math.\ Kyoto Univ. \textbf{19} (1979), 215--229.

%\bibitem{orlik}
%P.~Orlik, \emph{Seifert manifolds}, Lecture Notes in Mathematics 291,
%          Springer, Berlin, 1972.
%

\bibitem{rat}
J.~G. Ratcliffe.
\newblock \emph{Foundations of hyperbolic manifolds}, 
Springer-Verlag, New {Y}ork, 1994.

\bibitem{scot83}
P.~Scott, \emph{The geometries of 3-manifolds}, Bull. London Math. Soc.
  \textbf{15} (1983), 401--487.

\bibitem{thur}
W.~P. Thurston, \emph{The geometry and topology of 3-manifolds}, Princeton
  University Press.

%\bibitem{vais79}
%I.~Vaisman, \emph{Conformal foliations}, Kodai Math. J. \textbf{2} (1979),
%  26--37.
%
\bibitem{wol66}
J.~Wolf, \emph{Spaces of constant curvature}, McGraw Hill, New{Y}ork, 1966.

\bibitem{Woo86A}
J.C.~Wood, \emph{{H}armonic morphisms, foliations and {G}auss maps}, {C}omplex
  differential geometry and nonlinear partial differential equations (Y.T. Siu,
  ed.), Contemp.\ Math., vol.~49, Amer.\ Math.\ Soc., Providence, R.I., 1986,
  pp.~145--184.

\bibitem{Woo96}
J.C.~Wood, \emph{Harmonic maps and morphisms in $4$ dimensions},
Proceedings of the first Brazilian -- USA Workshop on Geometry,
Topology and Physics in Campinas, July, 1996 (to appear).

\end{thebibliography}
%\providecommand{\bysame}{\leavevmode\hbox to3em{\hrulefill}\thinspace}

\end{document}